\documentclass[12pt]{article}
\usepackage{amsfonts}

\title{Fano's inequality is also false for  three-dimensional quadric }
\author{Marat Gizatullin\thanks{Department of Mathematics,
 Technical University Federico Santa Mar{\'{\i}}a,
Avenida Espa\~na, No.\ 1680, Casilla 110-V, Valpara{\'{\i}}so,
Chile. {\emph{e-mail}} \textsf{mgizatul@mat.utfsm.cl}}}
\date{}

\newcommand{\PP}{\mathbb P}
\newcommand{\broken}{\dasharrow}
\newcommand{\mult}{\mathop{\rm mult}\nolimits}

\begin{document}
\bigskip
\par\noindent ABSTRACT. A birational transformation of the
tree-dimensional quadric is constructed.  The inequalities of Fano
are not fulfilled for this transformation.
\medskip
\maketitle
Acknowledgement. It is a pleasure to thank Miles Reid
who corrected the text of [3] . New variant of [3] is not yet
resubmitted to arXiv , but  Reid's lesson is used by the
preparation of the present paper.

\section*{Introduction}
It was shown in [3] that Fano's inequalities are false for
tree-dimensional projective space. The goal of this article is to
show the same thing for the case of three-dimensional quadric. Our
considerations are parallel to [3], and in a sense , an additional
text on the subject is not necessary , but if there exists an
uncountable set of papers based on the wrong inequality, then it
is worthy to publish at least two different texts explaining the
truth.
\medskip\par
Let $X$ be a non-singular threefold such that Pic$(X)=\mathbb{Z},$
and the anticanonical class \ $(-K_X) $\ is ample. If $(-K_X)=rH$
for a generator $H$ of the Picard group, then $X$ is said to be
Fano threefold of index $r$ ( and of the first kind ). According
to [1], Fano's inequality is  the statement:
 \begin{quote} \em
For any birational transformation,
 \[
 f\colon \ X \broken \ X
 \]
 defined by a linear system of  degree $d>1$
 (the degree is the
number defined by $f(mH)=dmH $ )
 either there exists a point $P\in X$ such that
 \[
 \mult_P(f(|H|))> 2d/r ,
 \]

 or there exists an irreducible curve $C\subset X$ such that
 \[
 \mult_C(f(|H|))> 2d/r .
 \]
\end{quote}
Let $Q\subset \mathbb{P}^4 $ be a smooth quadric hypersurface in
the four-dimensional projective space. $Q$ is a Fano threefold of
index 3. Therefore,  for this threefold, one can rewrite Fano's
inequality  as follows. \medskip \begin{quote} \em For any
birational transformation,
 \[
 f\colon \ Q \broken Q
 \]
 defined by five homogeneous polynomials of the same degree $d$ and
without a common two-dimensional set of zeros on $Q$,
 \[
 x_i'= f_i(x_0,x_1,x_2,x_3, x_4 ),\quad i=0,1,2,3,4,
 \]
 either there exists a point $P\in Q$ such that
 \[
 \mult_P(f_i)>2d/3
 \]
for every $i=0,1,2,3,4$, or there exists an irreducible curve
$C\subset Q$ such that
 \[
 \mult_C(f_i)> d/3
 \]
 for every $i$.
 \end{quote}

The goal of my article is to show that these inequalities do not
take place for a certain birational transformation  of degree
d=13. That is, I write down the formulas for a birational
transformation of degree 13 such that for the forms
$f_0,\dots,f_4$ defining the transformation, for any point $P\in
Q$ and for any irreducible curve $C\subset Q$ one can see that
 \[
 {\min}_i(\mult_P(f_i))\le8 \quad\hbox{and}\quad
{\min}_i(\mult_C(f_i))\le4.
 \]

\section*{Construction of the example}
Let us consider the homogeneous coordinates $x_0,x_1,x_2,x_3,x_4$
for $\PP^4$ as the normalized coefficients of a binary quartic
form $F(T_0,T_1)$,
 \[
 F(T_0,T_1)=x_{0}T_{0}^4+4x_{1}T_{0}^3T_{1}+6x_{2}T_{0}^2T_{1}^2+4x_{3}T_0T_{1}^3+x_4T_1^4,
 \]
 that is we consider
 $\PP^4$ as the projectivization of the vector space
of binary quartics. Let $I=I(x_0,x_1,x_2,x_3,x_4)$ be the
quadratic (or polar ) invariant of the binary quartic,
 \[
 I=x_{0}x_{4}-4x_{1}x_{3}+3x_{2}^2 ,
\]
$J=J(x_0,x_1,x_2,x_3,x_4)$ be the cubic  invariant (or the Hankel
determinant ) of the quartic,
\[
 J={\left |\begin{array}{ccc}
  x_{0} & x_{1} & x_{2} \\
  x_{1} & x_{2} & x_{3} \\
  x_{2} & x_{3} & x_{4}   \end{array}\right | } .
\]
$I$ and $J$ are ground ( or basic ) invariants of binary quartic,
therefore, if both of them vanish simultaneously on a quartic,
then the quartic is a null-form ( or unstable form ), that is it
has a triple linear factor.
\medskip \par
 We define our  quadric $Q$ with the help of the quadratic
  invariant,
  $$Q\quad :\qquad  I(x_0,x_1,x_2,x_3)=0.$$
\medskip \par
 We fix a parameter $t$ and consider five following forms of degree 13.
 \begin{eqnarray*}
 (f_t)_0&=&x_0J^4, \\
(f_t)_1&=&x_1J^4+tx_0^4J^3,\\
(f_t)_2&=&x_2J^4+2tx_1x_0^3J^3+t^2x_0^7J^2,\\
(f_t)_3&=&x_3J^4+3tx_2x_0^3J^3+3t^2x_1x_0^6J^2+t^3x_0^{10}J,\\
(f_t)_4&=&x_4J^4+4tx_3x_0^3J^3+6t^2x_2x_0^6J^2+4t^3x_1x_0^{9}J+t^4x_0^{13}.
 \end{eqnarray*}
These five forms define a one-parameter family of rational maps
 \[
 g_t\colon \PP^4\broken \PP^4.
 \]
 If $t=0$, all five forms have a common factor $J^4$. After cancelling
this, we see that $g_0$ is the identity transformation. For our
example, we need nonzero values of $t$. If $t$ is not zero, then
it is clear that the five forms have no nonconstant common factor.
Moreover,
 \[
 J\Big((f_t)_0,(f_t)_1,(f_t)_2,(f_t)_3,(f_t)_4\Big)=J(x_0,x_1,x_2,x_3,x_4)^{13}.
 \]
In fact, this identity expresses the invariance of the Hankel
determinant under triangular transformation of variables
$T_0,T_1$. Using the latter identity, it is not hard to see that
 \[
 (f_{-t})_i\Big((f_t)_0,(f_t)_1,(f_t)_2,(f_t)_3,(f_t)_4\Big)=x_iJ^{42},
 \]
that is,
 \[
 g_{(-t)}\circ g_t=\hbox{the identity transformation}.
 \]
 Thus $g_t$ is rationally invertible and is a Cremona transformation.
More generally,
 \[
 g_{s}\circ g_t=g_{s+t},
 \]
and we get a one-parameter group of Cremona transformations. These
transformations induce biregular automorphisms of an affine open
subset of the projective space, the complement of the
determinantal cubic hypersurface $J=0$. Indeed, the above formula
of the determinant transformation proves this (moreover, one can
see below the exact calculation of the fundamental points of such
a transformation).
\par Moreover
$$I\Big(x_0^{\prime},x_1^{\prime},x_2^{\prime},x_3^{\prime},x_4^{\prime}\Big)=I(x_0,x_1,x_2,x_3,x_4)\cdot J^{12},$$
this identity expresses the invariance of the quadratic invariant
under triangular transformation of variables $T_0,T_1$. Thus , the
quadric $Q$ is invariant by the transformations.
\medskip\par\noindent
{\textbf{Remark}}.

The formulas for $g_t$ can be generalized and some formulas for an
infinite dimensional family of automorphisms of the complement to
the determinantal hypersurface ( or induced automorphisms of the
complement of the determinantal hypersurface in the quadric ) can
be written. The generalization is similar to formulas   written
down on page 8 of the Max-Planck Institute preprint [2] and looks
as follows. If $\phi_m(x,y)$ is any binary form of degree $m$ ,
$\phi= \phi_m(x_0^3,J),$ then one can define the following
transformations .
 \begin{eqnarray*}
x_0^{\prime}&=&x_0J^{4m}, \\
x_1^{\prime}&=&x_1J^{4m}+x_0\phi J^{3m},\\
x_2^{\prime}&=&x_2J^{4m}+2x_1\phi J^{3m}+x_0\phi^2 J^{2m},\\
x_3^{\prime}&=&x_3J^{4m}+3x_2\phi J^{3m}+3x_1\phi^2
J^{2m}+x_0\phi^3
J^{m},\\
x_4^{\prime}&=&x_4J^{4m}+4x_3\phi J^{3m}+6x_2\phi^2
J^{2m}+4x_1\phi^3 J^{m}+x_0\phi^4.
 \end{eqnarray*}
 It is clear that
$$J\Big(x_0^{\prime},x_1^{\prime},x_2^{\prime},x_3^{\prime},x_4^{\prime}\Big)=J^{12m+1},$$
$$I\Big(x_0^{\prime},x_1^{\prime},x_2^{\prime},x_3^{\prime},x_4^{\prime}\Big)=I(x_0,x_1,x_2,x_3,x_4)\cdot J^{12m}.$$

\bigskip
\par Let us return to the family of transformations $g_t.$
We fix a nonzero value of the parameter $t$, for example, $t=1$,
and consider the corresponding Cremona transformation
 \begin{eqnarray*}
 x_0^{\prime}&=&x_0J^4, \\
x_1^{\prime}&=&x_1J^4+x_0^4J^3,\\
x_2^{\prime}&=&x_2J^4+2x_1x_0^3J^3+x_0^7J^2,\\
x_3^{\prime}&=&x_3J^4+3x_2x_0^3J^3+3x_1x_0^6J^2+x_0^{10}J,\\
x_4^{\prime}&=&x_4J^4+4x_3x_0^3J^3+6x_2x_0^6J^2+4x_1x_0^{9}J+x_0^{13}.
 \end{eqnarray*}
First of all, we find the points $P\in \PP^4$ where each form on
the right hand side has positive multiplicity (that is, the set of
all common zeros of these right hand sides, or, in other words,
the fundamental points of the transformation).

\par The first right hand side
vanishes if either $x_0=0$,  or $J=0$, or simultaneously $x_0=0,
J=0.$  If $x_0=0,$ but $J\neq 0,$ then
 using other four formulas, one  sees that for other four  coordinates of
a fundamental point, the equalities $x_1=x_2=x_3=x_4=0$ take
place. This case is not a point of $\PP^4$. \par The case $J=0,$
but $x_0\neq 0$ is also impossible for a fundamental point.
\par Thus , the  fundamental points consist of the
solutions of the following system of equations \[ x_0=0, \quad
J=0,\] in other words,
 the set of fundamental points is the hyperplane section (by the
hyperplane $H_0$ defined by $x_0=0 $) of the Hankel determinantal
threefold $D$.

The determinantal hypersurface $D$ has double points lying on the
twisted quartic $T$ parameterized by
 \[
 x_0=t_0^4,\quad x_1=t_0^3t_1,\quad x_2=t_0^2t_1^2, \quad
 x_3=t_0t_1^3,\quad x_4=t_1^4.
 \]
More precisely, the singular locus of the determinantal
hypersurface is $T$, and $\mult_PJ=2$ for every $P\in T$.
Therefore any  point not belonging to $T$ has multiplicity one on
the determinantal.
\par The next step. Let us consider the fundamental points of the
induced transformation $f|Q,$  that is  intersection
$$Q\cap D\cap H_0. $$
Intersection $Q \cap D$ consists of the quartic null-forms .
Binary quartic null-form $q$ has triple linear factor:
$q=b(T_0,T_1)\cdot a^3(T_0,T_1)$. If such unstable binary quartic
belongs to the hyperplane $H_0$, then the first coefficient $x_0$
of the quartic is zero,  $T_1$ is a divisor of the unstable
quartic . One can distinguish three cases:
 \begin{itemize}
\item either the triple factor $a^3$ is not proportional to $T_1^3,$
and $b$ is   proportional to $T_1$ ,\item  or $a=T_1,$ and $b$ is
not proportional to $T_1$,\item  or $q=T^4.$
\end{itemize}\medskip
\par Case 1: $q=T_1(uT_0+vT_1)^3,\quad u\neq 0.$
\par The set $B$ of these points is isomorphic to the affine line,
and it is clear that for this curve
$${\textnormal{mult}}_B(x_1J^4+x_0^4J^3)=4, $$
because  $B\neq T$, and  the points of $B$ are out of the
singularities of $D$.\medskip
\par Case 2: $q=(uT_0+vT_1)T_1^3,\quad u\neq 0.$
\par The set $C$ of these points is isomorphic to the affine line,
and it is clear that  for this curve
$${\textnormal{mult}}_C(x_3J^4+3x_2x_0^3J^3+3x_1x_0^6J^2+x_0^{10}J)=4. $$
\medskip
\par Case 3: $q=T_1^4 \quad $
$q$ is a double point of  hypersurface $D$, and
$${\textnormal{mult}}_q(x_4J^4+4x_3x_0^3J^3+6x_2x_0^6J^2+4x_1x_0^{9}J+x_0^{13})=8. $$

\par\noindent The Department of Mathematics,\par\noindent
 Technical University Federico Santa Mar{\'{\i}}a,\par\noindent
Avenida Espa\~na, No. 1680, Casilla 110-V,\par\noindent
Valpara{\'{\i}}so, Chile \medskip
\par\noindent {\emph{e-mail}} \textsf{mgizatul@mat.utfsm.cl}

  \end{document}